\documentclass[12pt, reqno]{amsart}
\usepackage{amssymb,amsthm,amsfonts,amsmath}

\usepackage{hyperref}
\usepackage{mathrsfs}

\newcommand\Y{\mathbb Y}
\newcommand\Z{\mathbb Z}
\newcommand\C{\mathbb C}
\newcommand\R{\mathbb R}

\newcommand\al{\alpha}
\newcommand\be{\beta}
\newcommand\la{\lambda}
\newcommand\om{\omega}
\newcommand\Om{\Omega}
\renewcommand\th{\theta}

\newcommand\ccdot{\,\cdot\,}

\begin{document}

\title[]
{The topological support of the z-measures\\ on the Thoma simplex}

\author{Grigori Olshanski}
\address{Institute for Information Transmission Problems, Bolshoy
Karetny 19,  Moscow 127051, Russia; Skolkovo Institute of science and Technology, Moscow, Russia; Department of Mathematics, National Research University Higher School of Economics}
\email{olsh2007@gmail.com}

\date{September 18, 2018}

\thanks{}

\begin{abstract}
The Thoma simplex $\Om$ is an infinite-dimensional space, a kind of dual object to the infinite symmetric group. The z-measures are  a family of probability measures on $\Om$ depending on three continuous parameters. One of them is the parameter of the Jack symmetric functions, and in the limit when it goes to $0$, the z-measures turn into the Poisson-Dirichlet distributions. The definition of the z-measures is somewhat implicit. We show that the topological support of any nondegenerate z-measure is the whole space $\Om$. 
\end{abstract}

\maketitle

\textbf{1.} The \emph{Thoma simplex}, denoted by $\Om$, is a closed subset of the double product space $[0,1]^\infty\times[0,1]^\infty$. The points of $\Om$ are pairs of sequences $\al=(\al_i)$, $\be=(\be_i)$ such that
$$
\al_1\ge\al_2\ge\dots\ge0, \quad \be_1\ge\be_2\ge\dots\ge0, \quad \sum_{i=1}^\infty \al_i+\sum_{i=1}^\infty\be_i\le1.
$$
The problem of harmonic analysis on the infinite symmetric group leads to a two-parameter family $\{M_{z,z'}\}$ of probability measures on $\Om$, called the \emph{z-measures}, see \cite{KOV1}, \cite{KOV2}, \cite{BO-1998}, \cite{BO-2000}. As first shown in \cite{Ker-2000} (see also \cite{BO-2005}, \cite{Ols-2010}), the definition of the z-measures can be extended by adding the third parameter $\th>0$, which leads to a three-parameter family $\{M_{z,z',\th}\}$ of z-measures. The values $\th=1$ and $\th=\frac12$ have special meaning: for $\th=1$ one recovers the two-parameter measures $M_{z,z'}$, and for $\th=\frac12$ one gets yet another family of measures related to harmonic analysis (see \cite{Str}).

\smallskip

\textbf{2.} The subset of $\Om$ formed by points with all $\be$-coordinates equal to $0$ can be identified with the \emph{Kingman simplex}, usually denoted by $\overline\nabla_\infty$.  On $\overline\nabla_\infty$, there is a well-known one-parameter family of probability measures --- the \emph{Poisson-Dirichlet distributions} $PD(\tau)$, $\tau>0$, and a more general family $\{PD(a,\tau)\}$ of Pitman's two-parameter Poisson-Dirichlet distributions, see \cite{Feng} (here we use a non-standard notation of the parameters  in order to avoid confusion with the notation in the previous paragraph). In the limit regime as parameters $z,z'$ go to infinity and parameter $\th$ goes to $0$ in such a way that $zz'\th\to\tau$, the z-measures converge to the Poisson-Dirichlet distribution $PD(\tau)$, see \cite[Section 1.2]{Ols-2010}. The two-parameter measures $PD(a,\tau)$ can also be obtained as a degeneration of z-measures, but in a more formal way, see \cite[Remark 9.12]{Ols-2010}. Although the z-measures and the Poisson-Dirichlet distributions have a number of similarities, the former are substantially more sophisticated objects than the latter. One of the differences between them is that the Poisson-Dirichlet distributions can be explicitly constructed using relatively simple probabilistic algorithms (see \cite{Feng}), while for the z-measures no similar construction is available.   

\smallskip

\textbf{3.} The z-measures can be defined in two ways: as limits of certain measures on Young diagrams of growing size and as the solution to a kind of moment problem (see \cite[Theorems 1.6 and 1.8]{BO-2005}). Let us briefly describe the second definition: the measure $M_{z,z',\th}$ is uniquely characterized by the equations
$$
\int_{\om\in\Om} P_\la(\om;\th)M_{z,z',\th}(d\om)=\varphi_{z,z',\th}(\la),
$$
where $\om$ is another notation for $(\al,\be)$; $\la$ ranges over the set of partitions $\Y$; $P_\la(\ccdot;\th)$ is the Jack symmetric function with parameter $\th$, realized as a function on $\Om$ via the specialization 
$$
p_1\mapsto1; \qquad p_k\mapsto \sum_{i=1}^\infty\al_i^k+(-\th)^{k-1}\sum_{i=1}^\infty\be_i^k, \quad k\ge2,
$$
where $p_k$ is the $k$th the power-sum symmetric function; finally, $\varphi_{z,z',\th}(\la)$ is a certain function on $\Y$.

The functions $P_\la(\om;\th)$ are nonnegative on $\Om$, so that $\varphi_{z,z',\th}(\la)\ge0$. We impose the condition of \emph{non-degeneracy}: our constraints on the parameters $(z,z')$ are chosen so that that $\varphi_{z,z',\th}(\la)>0$ for all $\la\in\Y$. See \cite[Proposition 1.2]{BO-2005} for more details. 

For more definiteness, we write down explicitly $\varphi_{z,z',\th}(\la)$ in the special case $\th=1$: then the Jack symmetric functions turn into the Schur functions and 
$$
\varphi_{z,z,1}(\la)=\dfrac{\prod_{(i,j)\in\la}(z+j-i)(z'+j-i)}{(zz')(zz'+1)\dots(zz'+n-1)}\cdot\frac{\dim\la}{n!}, 
$$
where the product in the numerator is over the boxes of the diagram of $\la$, $\dim\la$ is the number of standard tableaux of shape $\la$, $n=|\la|$ is the number of boxes in $\la$. In this case, the constraints on $(z,z')$ are the following: either $z'=\bar z\in\C\setminus\R$ or $m<z,z'<m+1$ for some $m\in\Z$. 

\smallskip

\textbf{4.} Our aim is to prove the following theorem

\smallskip

\noindent\textbf{Theorem.} \emph{Let $M_{z,z',\th}$ be an arbitrary nondegenerate z-measure. For any nonempty open subset $U\subset\Om$, its mass $M_{z,z',\th}(U)$ is strictly positive. In other words, the topological support of $M_{z,z',\th}$ is the whole space\/ $\Om$.}
\smallskip

Note that for the Poisson-Dirichlet distributions $PD(a,\tau)$, a similar claim also holds true: the topological support is the whole Kingman simplex. In fact, much more is known: the distribution of the first $n$ coordinates $\al_1,\dots,\al_n$ on $\overline\nabla_\infty$, $n=1,2,\dots$, is absolutely continuous with respect to Lebesgue measure, and there is  an expression for the density \cite[Theorem 3.6]{Feng}. For the three-parameter z-measures, results of such kind are unknown. 

We will show that the claim of the theorem is a direct consequence of the results of \cite{BO-2009}, \cite{Ols-2010}, and \cite{Kor} about a model of diffusion processes related to the z-measures. 

\begin{proof}
Let us abbreviate $M:=M_{z,z',\th}$. The construction of \cite{Ols-2010} assigns to $M$ 
a Feller semigroup $T(t):=T_{z,z',\th}(t)$ on $C(\Om)$ with a transition function $p(t;\om,\cdot)$, $t\ge0$ (see \cite[Theorem 9.6]{Ols-2010}). (Note that in the special case $\th=1$ the construction simplifies, it was described in \cite{BO-2009}.) 

Let $A$ denote the generator of the semigroup $T(t)$. As explained in \cite[Theorem 9.7]{Ols-2010} (with reference to  \cite[Theorem 8.1]{BO-2009}), for any point $\om\in\Om$ and its neighborhood $V$ one can exhibit a function $f\in C(\Om)$, contained in the domain of $A$ and such that $Af(\om)=0$, $\Vert f\Vert=f(\om)>0$, and the supremum of $f$ outside $V$ is strictly less than $f(\om)$. 

It is known (see \cite[ch. 4, Remark 2.10]{EK}) that this property of the generator implies the following:  for any open set $U\subset \Om$ and any point $\om\in U$, one has $1-p(t; \om, U)=o(t)$. In particular, $p(t;\om,U)>0$ for all $t$ small enough.

On the other hand, by virtue of  \cite[Theorem 17]{Kor}, the measure $p(t; \om, \ccdot)$  is absolutely continuous with respect to $M$, for any $\om\in\Om$ and $t>0$. We conclude that $M(U)>0$, as desired. 
\end{proof}

The result can be useful when working with the $L^2$ version of the Markov semigroup $T(t)$.

\end{document}